\documentclass[10pt]{article}

\usepackage{hhline}
\usepackage{amssymb, amsmath}
\usepackage{cite}
\usepackage{array}
\usepackage{theorem}
\usepackage{graphicx}
\theoremstyle{plain}

\newtheorem{pro}{Proposition}

\theoremstyle{remark}

\begin{document}

\renewcommand{\abstractname}{Abstract}
\renewcommand{\figurename}{Fig.}

\newcommand{\bfK}{{\bf K}}
\newcommand{\smj}{\suml _{j=1}^i}
\newcommand{\smn}{\suml _{j=1}^N}
\newcommand\os[1]{\trk{1.3}{#1}}
\newcommand{\bom}{{\bs\omega}}
\newcommand{\bgam}{{\bs\gamma}}
\newcommand{\bbt}{{\bs\beta}}
\newcommand{\bal}{{\bs\alpha }}
\newcommand{\bxi}{{\bs\xi}}
\newcommand{\blm}{{\bs\lm}}
\newcommand{\bOm}{{\bs\Omega}}
\newcommand{\bfJ}{{\bf J}}
\newcommand{\bfI}{{\bf I}}
\newcommand{\bfB}{{\bf B}}
\newcommand{\bfE}{{\bf E}}
\newcommand{\bfN}{{\bf N}}
\newcommand{\bfQ}{{\bf Q}}
\newcommand{\bfS}{{\bf S}}
\newcommand{\bfA}{{\bf A}}
\newcommand{\bfC}{{\bf C}}
\newcommand{\bfG}{{\bf G}}
\newcommand{\bfXi}{{\bf \Xi}}
\newcommand{\bfOm}{{\bf \Omega}}
\newcommand{\mcM}{\mathcal M}
\newcommand{\mcE}{\mathcal E}
\newcommand{\mcD}{\mathcal D}
\newcommand{\mcI}{\mathcal I}
\newcommand{\mcC}{\mathcal C}
\newcommand{\mcG}{\mathcal G}
\newcommand{\sgn}{\mathop{\rm sgn}\nolimits}
\newcommand{\ext}{\mathop{\rm ext}\nolimits}
\newcommand{\diag}{\mathop{\rm diag}\nolimits}
\newcommand{\rank}{\mathop{\rm rank}\nolimits}

\newcommand{\vs}[2]{#1_1,\ldots,#1_{#2}}

\newcommand\St{\mbox{St}}
\newcommand\ad{\mbox{ad}}
\newcommand\Ad{\mbox{Ad}}
\newcommand\goth{\mathfrak}
\newcommand\R{\mathbb R}
\newcommand\Z{\mathbb Z}
\newcommand\ind{\mathrm{ind}\,}
\newcommand\AV{}
\newcommand\RU{}
\newcommand\EN{}
\newcommand\NO{no.\,}
\newcommand\ndd{--}
\newcommand\bs{\boldsymbol}

\newcommand\proof{{\bf Proof.}\quad}

\newcommand\Ann{\mbox{Ann} }
\newcommand\spann{\mathrm{span} }
\newcommand\const{{\rm const}}

%%%%%%%%%%%%%%%%%%%%%%%%%

%\title{Topological monodromy in nonholonomic systems.}
%\author{A.\,V.\,Bolsinov, A.\,A.\,Kilin, A.\,O.\,Kazakov$^{1}$\\
%\normalsize \normalsize $^1$Institute of computer science, Izhevsk}
%\date{}
%\maketitle

\title{Topological monodromy as an obstruction to Hamiltonization of nonholonomic systems: \\ pro or contra?}
\author{A.\,V.\,Bolsinov$^{1}$, A.\,A.\,Kilin$^{2}$, A.\,O.\,Kazakov$^{2,3}$\\
\normalsize $^1$Department of Mathematical Sciences, Loughborough University, LE11 3TU, UK \\
\normalsize $^2$Institute of computer science,
\normalsize ul. Universitetskaya 1, Izhevsk, 426034, Russia \\
\normalsize $^3$ Research Institute of Applied Mathematics and Cybernetics, \\
\normalsize Nizhny Novgorod State University, pr. Gagarina 23, Nizhny Novgorod, 603950, Russia}
\date{}
\maketitle

\begin{abstract}%%
The phenomenon of a topological monodromy in integrable
Hamiltonian and nonholonomic systems is discussed.  An efficient method
for computing and visualizing the monodromy is developed.  The comparative
analysis of the topological monodromy is given for the rolling  ellipsoid
of revolution problem in two cases, namely, on a smooth and  on a rough
plane. The first of these systems is Hamiltonian, the second is
nonholonomic.  We show that,  from the viewpoint of monodromy, there is no
difference between the two systems,  and thus disprove the conjecture by
Cushman and Duistermaat stating that the topological monodromy gives a
topological obstruction for Hamiltonization of the rolling ellipsoid of
revolution on a rough plane.
\end{abstract}

{\bf Keywords.} {\em %ключевые слова
Topological monodromy,  integrable systems, nonholonomic
systems, Poincar\'e map, bifurcation analysis, focus-focus singularities.}

\section{Introduction}

The paper has been motivated by the following general question in classical mechanics:
how and to what extent does the dynamical behavior of nonholonomic systems differ from that of Hamiltonian ones?
This question is closely related to the Hamiltonization problem: is it possible to turn a given nonholonomic system into Hamiltonian by an appropriate choice of a Poisson
structure and change of time? This problem is quite nontrivial, discussed in many papers (see, e.g., \cite{BorMamRolling,
21,22,23,24,25,26,30,31,32,33,9*}) and has many aspects, one of which is finding topological obstructions to Hamiltonization of integrable nonholonomic systems.

Here by integrability we understand the existence of sufficiently many first integrals such that their  common regular levels are diffeomorphic to  two-dimensional tori (as in the
case of integrable Hamiltonian systems with two degrees of freedom).  The phase space of such a system is foliated into invariant 2-tori.  Speaking of topological
obstructions to Hamiltonization, we mean the following natural question: is it possible to find those properties of such a foliation which allow us to distinguish it from similar
foliations that appear in integrable Hamiltonian systems (the so-called Liouville foliations)?

Clearly, no such obstructions exist near a regular fiber. Moreover, it is well known that in the presence of an invariant measure the system (after an appropriate change of time)
admits a Hamiltonian representation (see \cite{kolm, kozlov1, BBM1}). However, topological obstructions may exist in a neighborhood of singular fibers. One of such obstructions is the so-called
topological monodromy of a foliation into invariant tori.  The difference between Hamiltonian and non-Hamiltonian monodromy was one of the main
issues studied in the famous paper by J.~Duistermaat and R.~Cushman \cite{CushDuis} where a detailed topological treatment of the monodromy in integrable nonholonomic systems was
given.  As a concrete example of a nonholonomic system, where the monodromy  is essentially non-Hamiltonian and  Hamiltonization is, therefore, impossible, the authors
suggest the problem of the rolling prolate ellipsoid of revolution on a rough plane (i.e., rolling without sliding).\footnote{Here is a citation from
\cite{CushDuis}: ``Because the monodromy going around this heteroclinic
cycle is the identity, the rolling prolate ellipsoid of revolution cannot
be made into a Hamiltonian system, even though it is time reversible and
energy conserving. This is an example where a global invariant (namely,
monodromy) has been used to show that a 4-dimensional conservative time
reversible system is not Hamiltonian''. Unfortunately, the paper does not contain any detailed explanations to this conclusion.}

However, it is well known that quite a similar problem in the case of a smooth plane (i.e., when the friction is zero) is Hamiltonian.  Thus, it would be very interesting to observe any difference in the dynamics of these two systems.
Since the monodromy is a rather rough topological characteristic,  the phenomenon should be easy to observe. Our preliminary considerations, however, did not
reveal any difference in the behavior of these systems and we decided to carry out a detailed  analysis of the topological monodromy for both of them.

The paper is organized as follows. In the next section we recall the notion of monodromy for integrable systems and discuss some of its properties in
the case of Hamiltonian systems. In particular, following \cite{CushDuis},
we make an emphasis on the difference between Hamiltonian and non-Hamiltonian cases.
Then we discuss one of possible methods for calculating monodromy in systems with rotational symmetry, which is based on analysis of some properties of the Poincar\'{e}
map for a specially chosen section. In Sections 3 and 4 we apply this method to study the monodromy in two integrable problems of a rolling prolate ellipsoid of revolution:
on a smooth plane (Hamiltonian case) and on a rough plane (nonholonomic case).

The main conclusion of our work is that from the viewpoint of monodromy these two systems behave absolutely in the same way. In particular, the monodromy does not give any
obstruction to Hamiltonization of this nonholonomic system. Moreover, our analysis shows, in fact, that the foliations into invariant tori in these two cases are
isomorphic.  However, this does not mean that the monodromy is useless for the Hamiltonization problem. On the contrary, it makes it possible to essentially
reduce the ``searching sector'' for a suitable Poisson structure. These conclusions are discussed in the closing section of the paper.

\section{Topological monodromy in integrable systems}\label{sec_bkk2}

The notion of a monodromy for integrable (Hamiltonian) systems was introduced by Duistermaat in~\cite{Duis} as one of obstructions to the existence of global action-angle variables.
Since this notion has a pure topological nature, i.e. it is completely defined by the properties of the foliation into invariant tori, we can easily extend it to the case of
nonholonomic integrable systems.

We recall the definition of monodromy in the case we are dealing with (some generalizations are discussed in~\cite{zhilin, efst}). Consider an integrable  system whose phase space is foliated into two-dimensional invariant submanifolds (tori). The singular fibers are ignored or just removed.
 Choose a particular torus $T_0$ and some deformation of it  $T_t$, $t\in [0,1]$, such that $T_1=T_0$. In other words,
 we consider a closed path in the space of parameters (i.e., values of the first integrals) that defines a deformation after which the torus returns to the initial position.

 Next we fix a pair of basis cycles
 $\lambda_0$, $\mu_0$ on the initial torus $T_0$ and, by changing them continuously in the process of deformation, we obtain a family of cycles  $\lambda_t$, $\mu_t$ forming a basis  on  $T_t$ for each fixed value of
$t\in [0,1]$. When the deformation is completed,
on the torus $T_1=T_0$ we obtain a pair of basis cycles $\lambda_1,
\mu_1$.   It is clear that if the deformation takes place inside a small neighborhood of  $T_0$, then the cycles so obtained are homologous to the initial cycles $\lambda_0, \mu_0$, i.e.   $\lambda_0$ and~$\lambda_1$ can be continuously deformed to each other inside $T_0$ (similarly for  $\mu_0$ and
$\mu_1$).  However if the family $T_t$ goes ``far'' from the initial torus $T_0$,  it may happen that new cycles $\lambda_1,
\mu_1$  are essentially different from $\lambda_0,
\mu_0$.   They nevertheless still form a basis and therefore, up to a homotopy, are related to the initial cycles by means of a certain integer unimodular matrix:
$$
\begin{pmatrix}
\lambda_1 \\[-2mm]
\mu_1
\end{pmatrix} =
\begin{pmatrix}
a & b \\[-2mm]  c & d
\end{pmatrix}
\begin{pmatrix}
\lambda_0 \\[-2mm]
\mu_0
\end{pmatrix}\!\!,  \qquad a,b,c,d \in \Z, \ ad-bc = 1.
$$

That is exactly what is called the {\it monodromy matrix} corresponding to the deformation $T_t$, $T_0=T_1$. If it is different from the identity matrix we say that the monodromy is non-trivial.

Let us make some general comments about the monodromy which clarify its nature.

If we consider the foliation of the phase space $\mathcal M^4$ into invariant manifolds\footnote{This construction does not change if we consider a dynamical system  on a five-dimensional space $\mathcal M^5$ which admits three integrals $H, F_1, F_2$.  The rolling ellipsoid on a plane is a system of this kind.}, related to two integrals  $H$ and~$F$,
then it is convenient   to consider the integral map
$\Phi=(H, F)\colon  \mathcal M \to \R^2$, its image $\Phi(\mathcal M)$
and the bifurcation diagram $\Sigma \subset \Phi(\mathcal M) \subset \R^2$.
Then  choosing an initial torus  $T_0$ is equivalent to choosing a non-singular  (that is lying outside of $\Sigma$)  point $a\in \Phi(\mathcal M)$.   The torus $T_0$ itself is the preimage of $a$.  The deformation of the torus is defined by choosing a closed curve $\gamma (t)$ in the image of the integral map which does not intersect the bifurcation diagram (here we, or course, assume that
$\gamma(0)=\gamma(1)=a$).  The curve $\gamma$ defines a deformation of the torus
$T_t=\Phi^{-1}(\gamma(t))$  and, consequently, the monodromy.

If the curve $\gamma$  in the image of the momentum map is continuously deformed
in such a way that the deformation does not touch the bifurcation diagram, then the monodromy won't change. In particular, a non-trivial monodromy may appear for non-contractible loops $\gamma$ only.
Such non-contractible curves do not always exist, but very often they do, in particular, if
the bifurcation diagram contains isolated singular points. In this case, as a non-contactible loop one can take a small circle around the point. That is exactly the situation we are interested in. Let us discuss it in more detail.

\begin{figure}[!ht]
\centering
\includegraphics{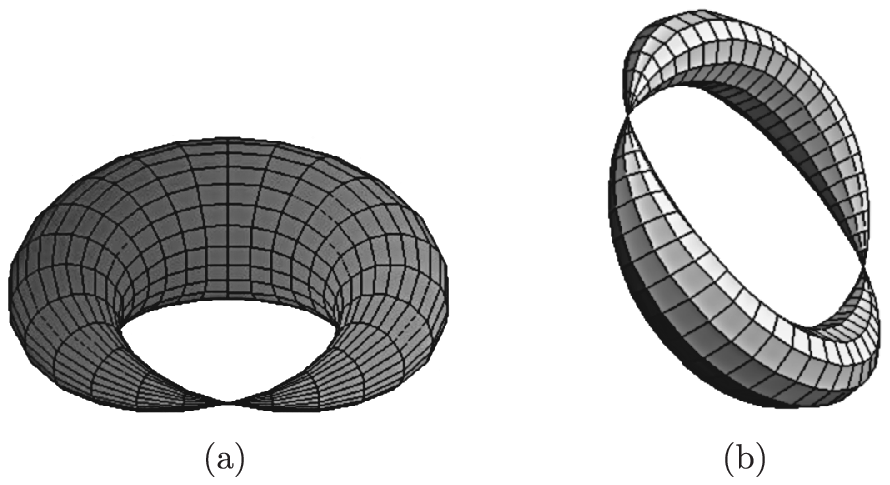}
\caption{~}\label{fig_N}
\end{figure}

In integrable systems one often deals with the situation when a singular integral manifold is a torus with a pinch point  or several pinch points (see Fig.~\ref{fig_N}).
The singular points in this case are equilibria of the system and are of focus-focus type. The topology of such singularities has been systematically studied in a range of papers~\cite{Matvfocus,Zungfocus,CushDuis,lerman,izosimov}, where the reader will find a detailed explanation of necessary definitions and results (see also \cite{bols_osh, bols_fom, cushbates, efst}). We restrict ourselves with a list of main properties of focus singularities which will be essential to understand our construction below.

\begin{itemize}
\item  Focus singular fibers are isolated in the sense that all of their neighboring fibers in the phase space~$\mcM^4$  are non-singular, i.e. are diffeomorphic to tori. On the bifurcation diagram of the integral map $\Phi$ they occur as isolated points. Typical examples are the Lagrange top, Clebsch case, spherical pendulum and ``champagne bottle'' system (see, e.g., \cite{cushbates, bols_fom, efst}).

\item  Going around such a singular point in the image of the integral map defines a non-trivial deformation of the torus,  and one can ask the question about the monodromy. Each focus point located of the singular fiber gives the same contribution to the monodromy, namely the matrix of the form
$$
\begin{pmatrix}
1 & 1 \\[-2mm] 0 & 1
\end{pmatrix}\!\!.
$$
In the Hamiltonian case (see  \cite{Matvfocus}, \cite{Zungfocus}) these matrices
have just to be multiplied so that the total monodromy matrix takes the form
$$
\begin{pmatrix}
1 & 1 \\[-2mm] 0 & 1
\end{pmatrix}^k =
\begin{pmatrix}
1 & k \\[-2mm] 0 & 1
\end{pmatrix}\!\!,
$$
where $k$ is the number of singular points on the fiber (that is the number of pinch points).
In the non-Hamiltonian case, the situation might be different  (see
\cite{CushDuis}). Namely, each point gives the same contribution as above, but as a factor in the total product one should take either this matrix itself or its inverse. In other words, these ``partial monodromy'' matrices may, as a result,  compensate each other. For example, in the case of a double pinched torus
(Fig.~\ref{fig_N}b) the total monodromy may appear to be trivial, which is impossible in the Hamiltonian case.

\item Focus singularities are stable under (integrable) perturbations of a system.
In particular, if a system depends on a parameter $J$, then the focus singularities will be preserved (``survive'') for all sufficiently small values of the parameter.
In particular, if instead of a four-dimensional phase space we consider a family of four-dimensional phase spaces parametrized by a Casimir function (or just by a certain integral in the nonholonomic case), then we obtain a family of focus singularities. In this case, it is convenient to consider the ``three-dimensional'' integral map by adding the additional parameter $J$ to the integrals
$H$, $F$.  On the bifurcation diagram of such a map,  focus singularities will occur as a curve. The monodromy makes sense in this case too, and  it is important that in the process of deformation, when one goes around the curve, it is allowed to vary the values of all the integrals including $J$ (i.e. $J$ need not to remain constant under deformation).

\end{itemize}

A typical example when a focus singularity necessarily appears is a natural system with a Hamiltonian  $H= K + V$ which is rotationally invariant.
Every non-degenerate equilibrium point that corresponds to a local maximum of the potential $V$ is a focus singular point. A similar statements holds true for nonholonomic systems too, i.e. a rotational symmetry naturally leads to the appearance of singularities of focus type (in the non-Hamiltonian case, such singularities have been studied in detail in~\cite{CushDuis}).

If we want the singular fiber to contain not one, but two focus points, we need to consider a dynamical system possessing an additional
$\Z_2$-symmetry.  Such a situation appears in the rolling problem for a prolate ellipsoid of revolution on the plane (both smooth and rough, i.e. with friction and without). When the ellipsoid takes the vertical position, we get an unstable equilibrium point. Since the system is rotationally invariant, this point is of focus type. Since the highest and lowest points (top and bottom) of the ellipsoid are symmetric, we get two distinct focus points.
Moreover, they belong to the same integral surface. This follows from the fact that if we slightly push the ellipsoid (standing at the vertical position) then first  it ``falls down'' but then returns to the vertical position again in such a way that the ``top'' and ``bottom'' interchange.

If the ellipsoid rolls on a smooth plane (i.e. without friction), then this dynamical system is known to be Hamiltonian, and the monodromy is given by the matrix
$A_{\mathrm{monodr}}=\begin{pmatrix}
1 & 2 \\[-2mm]
0 & 1
\end{pmatrix}\!.
$

However, this is still an open question whether or not the similar system on a rough plane (i.e. rolling without sliding) is Hamiltonian, and therefore according to the general non-Hamiltonian monodromy theorem  \cite{CushDuis},
a priori there are two possibilities:
\begin{equation}
\label{twovariants}
A_{\mathrm{monodr}}=\begin{pmatrix}
1 & 2 \\[-2mm]
0 & 1
\end{pmatrix}
\quad\mbox{
or}\quad
A_{\mathrm{monodr}}=\begin{pmatrix}
1 & 0 \\[-2mm]
0 & 1
\end{pmatrix}\!\!.
\end{equation}

If the second one takes place\footnote{That is exactly what is stated in~\cite{CushDuis}.}, then the Hamiltonization of this nonholonomic systems is surely impossible. But which of these possibilities takes place in reality?

The answer to this question will be given below and we shall see that the question itself is not well posed: the result depends on the way we go around this singularity. The point is that in the rolling ellipsoid problem, the space in which the system is naturally defined is five-dimensional and the system possesses 3 independent integrals. Thus the situation is more complicated than that in the model example when one goes around an isolated point on the plane. Now the ``going around process'' should be carried out in dimension 3 and, as we shall see below, two essentially different scenarios are possible.

\subsection{How are we going to compute the monodromy?}\label{subsec2_1_}

Usual methods for computing the monodromy are based on a combination of analytic and topological arguments. First one finds explicit analytic formulas for the integrals and after that these formulas are analysed by using some topological tools.

In the problem of a rolling ellipsoid on a rough plane, such a method does not work as explicit formulas for the integrals are unknown.  Instead we suggest to {\it visualize} the monodromy,  rather than to {\it compute} it,  by using the dynamics and numerical integration.
The main idea is that the monodromy can be reconstructed from some very natural dynamical properties of the system (a similar approach has been developed in
~\cite{ddsz}).  Notice that the construction presented below can be modified for a much wider class of dynamical systems including non-integrable ones.
The phenomenon we are going to observe is not related directly to integrability.  In fact, this is a certain property of the Poincar\'e map which is quite ``rigid'' and hence ``survives''  under small perturbations of the system (including, of course, those which are non-integrable). In the integrable case in question, this property can be interpreted in terms of the topological monodromy, and this interpretation is one of the key points of our paper.

Thus, we want to study the monodromy that is related to a ``walk'' around
an isolated singular point in the image of the integral map $\Phi=(H,
F)\colon  \mathcal M^4 \to \R^2$  (or,  which is the same, to a
``walk'' around an isolated singular fiber in the phase space). The main
example we are interested in is a singular fiber of focus type (see
Fig.~\ref{fig_N}), but the construction below can work in a more general
situation.

Let  $\gamma$  be a closed path around a singular point in the image of the integral map $\Phi$.  As  $\gamma$ we can take a circle of sufficiently small radius with the  angle
$\alpha\in [0, 2\pi]$ as a parameter on it. Consider the preimage of $\gamma$ under the integral map $\Phi$. Then each point of the curve $\gamma$ corresponds to an invariant torus and, hence,  the preimage of $\gamma$, as a whole, is a three-dimensional manifold  $Q^3$.
From the topological viewpoint $Q^3$ is a $T^2$-fiber bundle over the circle $\gamma$.

The topology of $Q^3$ is easy to describe just by using the monodromy
$A_{\mathrm {monodr}}$.  One should first take the direct product  $T^2
\times [0,2\pi]$,  and then glue the top and bottom of this cylinder, i.e.  identify the tori
 $T^2\times \{0\}$ and~$T^2\times \{ 2\pi \}$ by means of the ``linear'' map which is given in standard angle coordinates by the matrix
  $A_{\mathrm {monodr}}$.  Our goal is to reconstruct this matrix.

  The solution of this problem  can be essentially simplified in the case when the system admits an additional $SO(2)$-symmetry.
  Focus singularities always satisfy this property (see \cite{Zungfocus}), and in our problems about rolling ellipsoids the rotational symmetry has an obvious geometric reason,  as we consider an ellipsoid of {\it revolution}. We will use the presence of the symmetry vector field  (generator of the $SO(2)$-action) whose trajectories are all closed and lying on integral surfaces of the dynamical system.

  This leads immediately to a special form of the monodromy matrix. Indeed, according to the general scheme we need to choose a pair of basis cycles $\lambda_0, \mu_0$ on an invariant torus and then look after their evolution while ``walking around'' the singular fiber.
  In the presence of a $SO(2)$-symmetry vector field,  as one of these basis cycles, say $\mu_0$, we can choose  a trajectory of this symmetry field.  Moreover,
  this can be done simultaneously for all the tori, i.e. we may assume  that $\mu_\alpha$, $\alpha\in[0,2\pi]$ is always a trajectory of the symmetry field.  Thus, after completing the ``walking around'' process we get $\mu_{2\pi}=\mu_0$,  which means that the monodromy matrix takes the form
  $$
A_{\mathrm {monodr}}=\begin{pmatrix}
1 & k \\[-2mm] 0 & 1
\end{pmatrix}
$$
for an arbitrary choice of the other additional cycle  $\lambda$,  so the question is just to find one single number $k\in \Z$.

\begin{figure}[!ht]
\centering\includegraphics{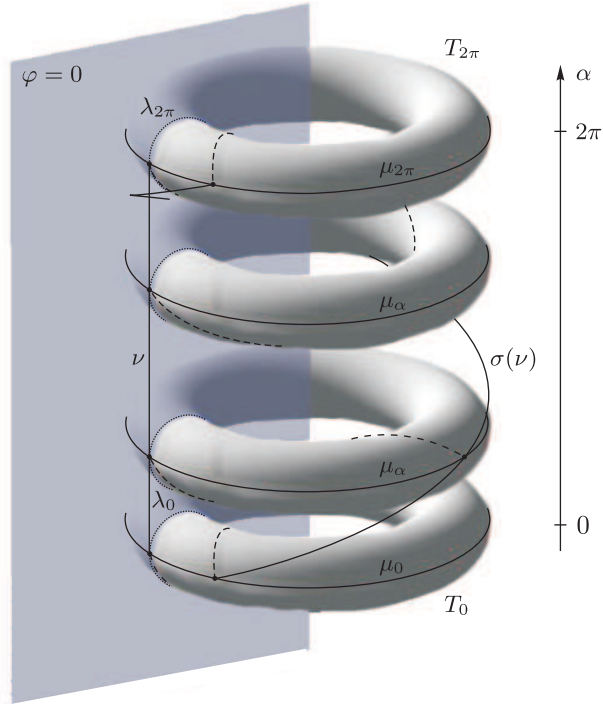}
\caption{Manifold $Q_3$ with basis cycles $\mu_{\alpha }$, $\nu$ and~$\sigma(\nu)$ on it.}
\label{fig3_}
\end{figure}

To that end, we use the fact that apart from the structure of a fibration into invariant tori, on $Q^3$ there is  an additional structure, namely
the initial dynamical system.  We illustrate our idea with a picture (see Fig.~\ref{fig3_}).  On this figure, the direct product  $T^2 \times
[0,2\pi]$ is shown with basis cycles  $\lambda_i$, $\mu_i$,
$i=0,2\pi$, indicated on the top and bottom bases.  The dotted lines illustrate phase trajectories of the system.  In order to obtain $Q^3$, one  should glue the top and bottom bases of this cylinder by means of a linear map which is given in the indicated bases by the formula
$$
\begin{pmatrix}
\lambda_{2\pi} \\[-2mm]
\mu_{2\pi}
\end{pmatrix} =
\begin{pmatrix}
1 & k \\[-2mm]  0 & 1
\end{pmatrix}
\begin{pmatrix}
\lambda_0 \\[-2mm]
\mu_0
\end{pmatrix}\!\!,
$$
where $k\in \Z$ (this number is to be found).  Here  $\mu$ is the uniquely defined cycle which, under the gluing operation,  is mapped to itself. In particular, if we consider the cylinder connecting the cycles  $\mu_0$
and~$\mu_{2\pi}$ (in other words,  we consider the whole family of cycles
$\mu_\alpha$, $\alpha\in[0,2\pi]$),  then after gluing the bases
$T^2\times \{0\}$ and~$T^2 \times \{2\pi\}$,  the cycles $\mu_0$
and~$\mu_{2\pi}$ will be identified,
and as a result this cylinder becomes a two-dimensional torus
$T_{\mathrm{transv}}$ lying inside  $Q^3$ and composed of  $\mu$-cycles.

This torus is not invariant. On the contrary, it is transversal to the dynamical system and, therefore, we can consider it as a global Poincar\'e section for the flow on
$Q^3$ (cf.~\cite{bdw}).

Now let us look at the topology of $Q^3$ from the other side. It is not hard to see that by cutting  $Q^3$ along the torus $T_{\mathrm{transv}}$, we again obtain the direct product
 $T^2\times I $  (where $I$ is some interval), but now the role of  $T^2$  is played not by an invariant torus but by the transversal torus $T_{\mathrm{transv}}$.   From a purely
topological viewpoint, this observation is equivalent to saying that $Q^3$ possesses two different structures of a $T^2$-fiber bundle over a circle.  In the first case, the fibers are invariant tori,  in the second one, fibers are tori transversal to the flow.

Consider the second (transversal) representation.  What is the gluing map in this case?
For the transversal torus the answer is evident: this is just the Poincar\'e map defined by the flow:
\begin{equation}
\sigma \colon  T_{\mathrm{transv}}  \to T_{\mathrm{transv}}.
\label{Poincaremap}
\end{equation}

This map also transforms basis cycles
on the Poincar\'e section  $T_{\mathrm{transv}}$.
 It remains to notice that the matrix of this transformation is the same as the monodromy matrix we are interested in. This fact can be ``seen'' directly
 (see~Fig.~\ref{fig3_});  it also follows from purely topological arguments (for instance, from computing the fundamental group $\pi_1(Q^3)$ by using two different representations of $Q^3$ as $T^2$-fiber bundles).

The main conclusion from this reasoning is as follows:

{\it In the presence of a rotational symmetry, we can find a global transversal Poincare surface of section $T_{\mathrm{transv}}\subset Q^3$  which is diffeomorphic to a torus.
The desired monodromy matrix  coincides with the  transformation matrix of basis cycles defined by the corresponding Poincar\'e map~\eqref{Poincaremap}}.

In other words, the monodromy can be found by analysing the properties of the Poincar\'e map defined on an appropriately chosen surface of section.

Since in our case the monodromy matrix has only one essential entry  $k\in \Z$, our problem can be solved by ``visualizing'' this number $k$ in the following way.
We choose two basis cycles on
$T_{\mathrm{transv}}$. One of them is a closed trajectory of the symmetry field
which has been earlier denoted by $\mu$.  Let $\nu$ be an arbitrary additional cycle.
Apply the Pincare map to both cycles  $\mu$ и~$\nu$.  The cycle $\mu$,  being a trajectory of the symmetry field, does not change, whereas
 $\nu$ is mapped to~$\nu + k
\mu$, i.e. its image $\sigma(\nu)$ will pass $k$ times along
 $\mu$.  This number of passages $k$ can be clearly seen if, on the torus $T_{\mathrm{transv}}$, we draw the image $\sigma(\nu)$ of the cycle $\nu$ under the Poincar\'e map
  \eqref{Poincaremap}.

That is exactly the method of computing (visualizing) the monodromy that we realize below in the rolling problem for an ellipsoid of revolution.
To carry this program out we only need to choose an appropriate torus $T_{\mathrm{transv}}$ transversal to the flow and then to evaluate numerically the action
of the Poincar\'e map on the basis cycles. This approach is accomplished in the two next sections.

\section{Rolling of an ellipsoid of revolution on a smooth plane}
\label{sec_bkk1}

\subsection{Equations of motion and first integrals}

Consider  a dynamically and geometrically axisymmetric ellipsoid rolling
on a smooth plane in a gravitational field. We assume that its center of
mass $O$ coincides with the geometrical center and choose a moving
coordinate system $O{\bs e}_1{\bs e}_2{\bs e}_3$ whose axes coincide with
the principal axes of inertia of the ellipsoid, and let $m$ denote its
mass and $\bfI=\diag (I_1,I_1,I_3)$~the central tensor of inertia. Then in
the chosen coordinate system, the equation of the ellipsoid takes the form
$({\bs r}, \bfB^{-1}{\bs r})=1$, where $\bfB=\diag(b_1^2, b_1^2, b_3^2)$,
and~$b_1$ and~$b_3$~are the principal semi-axes of the ellipsoid. Here and
in the sequel all vectors and tensors are assumed to be given in the
moving axes ${\bs e}_1$, ${\bs e}_2$, ${\bs e}_3$.\looseness=-1

\begin{figure}[!ht]
\centering\includegraphics{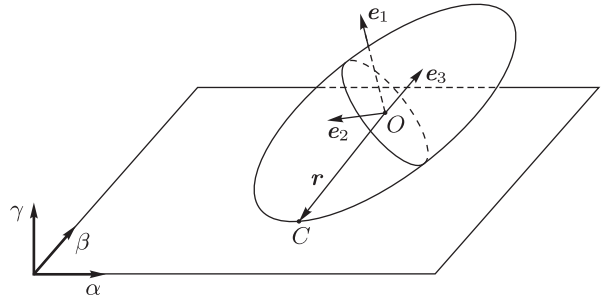}
\caption{An ellipsoid on a smooth plane.}\label{fig_1}
\end{figure}

As is well known~\cite{dtt}, for a body rolling on a smooth plane in a gravitational field the equations governing the evolution of the body's angular momentum ${\bs M}$ relative to the point of
contact  and the normal vector to the plane $\bgam$ can be written in closed Hamiltonian form as
\begin{equation}
\label{eq1_1}
\dot{\bs M}={\bs M}\times\dfrac{\partial H}{\partial {\bs M}}+\bgam\times\dfrac{\partial H}{\partial \bgam},\quad
\dot\bgam=\bgam\times\dfrac{\partial H}{\partial \bgam},
\end{equation}
where
\begin{equation}
\label{eq1_2}
H=\dfrac12 (\bfI\bfA{\bs M},\bfA{\bs M})+\dfrac12 ({\bs a}, \bfA{\bs M})^2-mg({\bs r},\bgam),
\end{equation}
${\bs a}=\bgam\times{\bs r}$, $\bfA=(\bfI+m{\bs a}\otimes {\bs a})^{-1}$, $g$~is the free-fall acceleration and the vector~${\bs r}$ directed from the point of contact to the center of mass (see Fig.~\ref{fig_1})
is related with~$\bgam$ by
\begin{equation}
\label{eq1_3}
{\bs r}=\dfrac{-\bfB\bgam}{\sqrt{(\bgam,\bfB\bgam)}}.
\end{equation}
The evolution of the remaining two unit vectors of the fixed coordinate system $\bal$ and $\bbt$ is given by quadratures
$$
\dot\bal=\bal\times\dfrac{\partial H}{\partial \bgam},\quad
\dot\bbt=\bbt\times\dfrac{\partial H}{\partial \bgam}.
$$

The Poisson bracket of Eqs.~\eqref{eq1_1} is defined by the algebra $e(3)$
$$
\{M_i,M_j\}=-\varepsilon_{ijk} M_k,\quad
\{M_i,\gamma _j\}=-e_{ijk}\gamma _k,\quad
\{\gamma _i,\gamma _j\}=0
$$
and possesses two Casimir functions
\begin{equation}
\label{eq1_4}
F_1=({\bs M},\bgam), \quad
F_2=\bgam^2,
\end{equation}
which are the integrals of Eqs.~\eqref{eq1_1}.
The integral~$F_2$ is the square of the unit normal vector $\bgam$ and is always equal to 1.
Hence,  the system~\eqref{eq1_1} defines the flow in the five-dimensional phase space $\mcM^5=\{{\bs M},\bgam\colon \bgam^2=1\}$.

It is well known that Eqs.~\eqref{eq1_1} admit another integral of motion (Lagrange integral) related to the rotational symmetry
\begin{equation}
\label{eq1_5}
F_3=M_3,
\end{equation}
and, therefore, are integrable. Thus, the phase space $\mcM^5$
is foliated (almost everywhere) into two-dimensional invariant manifolds, which are parameterized by the values of the integrals $H=h$, $F_1=p_{\psi}$,
$F_3=p_{\varphi}$. Here and throughout the paper $H$, $F_1$ and~$F_3$ are regarded as  functions of the phase variables $({\bs M}, \bgam)$, and $h$, $p_{\varphi}$ and
$p_{\psi}$~are constant values of these functions, which remain the same during the motion along a specific phase trajectory.

\subsection{Bifurcation analysis}\label{subsec_bkk1_2}

Two-dimensional (at fixed values $F_1=p_{\psi}$) bifurcation diagrams for an ellipsoid of revolution with the center of mass displaced along the symmetry axis were constructed by
M.\,Yu.\,Ivochkin in~\cite{ivochkin}.  Though the three-dimensional bifurcation diagram that we need for our analysis can be easily derived from~\cite{ivochkin},  we construct it below by using a different method as for our purposes it is important to draw an analogy to the
case that will be examined in Section~\ref{sec4} in a similar way.

Consider a three-dimensional integral map $\Phi=(p_\varphi,p_\psi,h)\colon
\mcM^5\,{\to}\, {\mathbb R}^3$ and the corresponding three-dimensional bifurcation
diagram~$\Sigma$, which is the image of the critical set of $\Phi$
\[
S =\{x\colon\rank d\Phi(x)\,{<}\,3\}
\]
and consists of two subsets $S_2 =\{x\colon\rank
d\Phi(x)=2\}$ and~$S_1 =\{x\colon\rank d\Phi(x)=1\}$.
The set~$S_2 $ is a two-parameter family of closed curves in~$\mcM^5$, which (in the typical case) are
periodic solutions of~\eqref{eq1_1}. The image of this set in the space of first integrals is a bifurcation surface. The set~$S_1 $ consists of two one-parameter
families of equilibrium points in~$\mcM^5$, and the corresponding image in the~space of first integrals consists of two curves.

\begin{pro}
For the case of a dynamically symmetric ellipsoid of revolution rolling on a smooth plane, the bifurcation diagram in the space of first integrals $(p_{\varphi}, p_{\psi}, h)$ consists of

1. the surface of regular precessions given by the equation
\begin{equation}
\label{eq1_6}
h=\dfrac{(p_\psi-p_\varphi\gamma _3)^2}{2I_1(1-\gamma ^{2}_{3})}+
\dfrac{p^{2}_{\varphi}}{2I_3}+mg\sqrt{b^{2}_{1}-(b^{2}_{1}-b^{2}_{3})\gamma ^{2}_{3}},
\end{equation}
where $\gamma _3$ is the solution of the equation
\[
\dfrac{(p_\psi-p_\varphi\gamma _3)(p_\psi\gamma _3-p_\varphi)}{I_1(1-\gamma ^{2}_{3})^2}+
\dfrac{mg(b^{2}_{3}-b^{2}_{1})}{\sqrt{b^{2}_{1}-(b^{2}_{1}-b^{2}_{3})\gamma ^{2}_{3}}}\,
\gamma _3=0,
\]

2. two families of (relative) equilibrium points determined by the curves
\begin{equation}
\label{eq1_7}
h=mg b_3 +\dfrac{p_{\varphi}^2}{2I_3},\quad
p_{\psi}=\pm p_{\varphi}.
\end{equation}
\end{pro}

\proof Choose the variables ${\bs x}=(\gamma _3, \dot\gamma _3, \varphi,p_{\varphi},p_{\psi})$ as local coordinates on $\mcM^5$, where $\varphi$~is the angle of self-rotation,
$p_{\varphi}=M_3$ and~$p_{\psi}=({\bs M},\bgam)$~are the momenta canonically conjugate
to the angles of intrinsic rotation and precession (which are integrals of motion in this case), $\gamma _3$~is the projection of $\bgam$ onto the symmetry axis of the
body, which is related to the nutation angle by $\gamma _3=\cos\theta $. The chosen coordinates are defined everywhere on $\mcM^5$ except for those points where the equality
$\gamma _3=\pm1$ holds; we shall consider these points separately.

We first consider the case where the rank of the map~$\Phi$ drops by 1
(the critical set~$S_2 $). For this purpose, we write in the chosen variables
the corresponding Jacobian
\begin{equation}
\label{eq1_8}
\dfrac{\partial (H,F_1,F_3)}{\partial {\bs x}}=\begin{pmatrix}
\frac{\partial H}{\partial \gamma _3} & \frac{\partial H}{\partial \dot\gamma _3} & \frac{\partial H}{\partial \varphi} & \frac{\partial H}{\partial p_{\varphi}} & \frac{\partial H}{\partial p_{\psi}}  \\
0 & 0 & 0 & 0 & 1\\
0 & 0 & 0 & 1 & 0
\end{pmatrix}\!\!,
\end{equation}
where the Hamiltonian has the form
\begin{equation}
\label{eq1_9}
H=\dfrac12 \left(\dfrac{m(b_1^2-b_3^2)\gamma _3^2}{b_1^2-(b_1^2-b_3^2)\gamma _3^2}+\dfrac{I_1}{1-\gamma _3^2}\right)\dot\gamma _3^2
+\dfrac{(p_\psi-p_\varphi\gamma _3)^2}{2I_1(1-\gamma _3^2)}+\dfrac{p_\varphi^2}{2I_3} +mg\sqrt{b_1^2-(b_1^2-b_3^2)\gamma _3^2}.
\end{equation}
Since $\frac{\partial H}{\partial \varphi}=0$, the condition for the rank to fall with $\gamma _3\ne\pm1$ can be represented as two algebraic equations
\begin{equation}
\label{eq1_10}
\dfrac{\partial H}{\partial \gamma _3}=0,\quad
\dfrac{\partial H}{\partial \dot\gamma _3}=0,
\end{equation}
which define in~$\mcM^5$ a two-parameter (with parameters~$p_\varphi$ and
$p_\psi$) family of closed curves forming~$S_2 $. These curves are periodic solutions of the system~\eqref{eq1_1} and
are called regular precessions (since they correspond to rotations of the body with a constant angle of inclination of the axis of rotation relative to the vertical). Substituting
the solution of Eqs.~\eqref{eq1_10} into the Hamiltonian~\eqref{eq1_9}, we obtain the bifurcation surface of regular precessions~\eqref{eq1_6} in the space of first integrals.

Now consider a neighborhood of the point $\gamma _3=\pm 1$.
As local coordinates near $\gamma _3=\pm 1$ we choose the variables ${\bs y}=(\gamma _1,\gamma _2,M_1,M_2,M_3)$. The Jacobian~\eqref{eq1_8} in these coordinates
takes the form
\begin{equation}
\label{eq1_11}
\left.\dfrac{\partial (H,F_1,F_3)}{\partial {\bs y}}\right|_{\gamma _3=\pm 1}=\begin{pmatrix}
\left.\frac{\partial H}{\partial \gamma _1}\right|_{\gamma _3=\pm 1} & \left.\frac{\partial H}{\partial \gamma _2}\right|_{\gamma _3=\pm 1} &
I^{-1}_{1}M_1 & I^{-1}_{1}M_2 & I^{-1}_{1}M_3 \\
M_1 & M_2 & 0 & 0 & \pm1\\
0 & 0 & 0 & 0 & 1
\end{pmatrix}\!\!.
\end{equation}
It is straightforward to notice that the rank of matrix~\eqref{eq1_11} drops (and by 2 at once) for $M_1=M_2=0$. The corresponding critical set is a one-parameter
family of relative equilibria of the system~\eqref{eq1_1} and can be written as $S_1
=\{{\bs M}=(0,0,p_\varphi),\bgam=(0,0,\pm1)\colon p_\varphi={(-\infty,+\infty)}\}$.
Physically these equilibria are vertical rotations of the ellipsoid about its axis of symmetry.
The image of this set in the space of first integrals consists of two bifurcation curves~\eqref{eq1_7}.$\blacksquare$

\begin{figure}[!ht]
\centering
\includegraphics{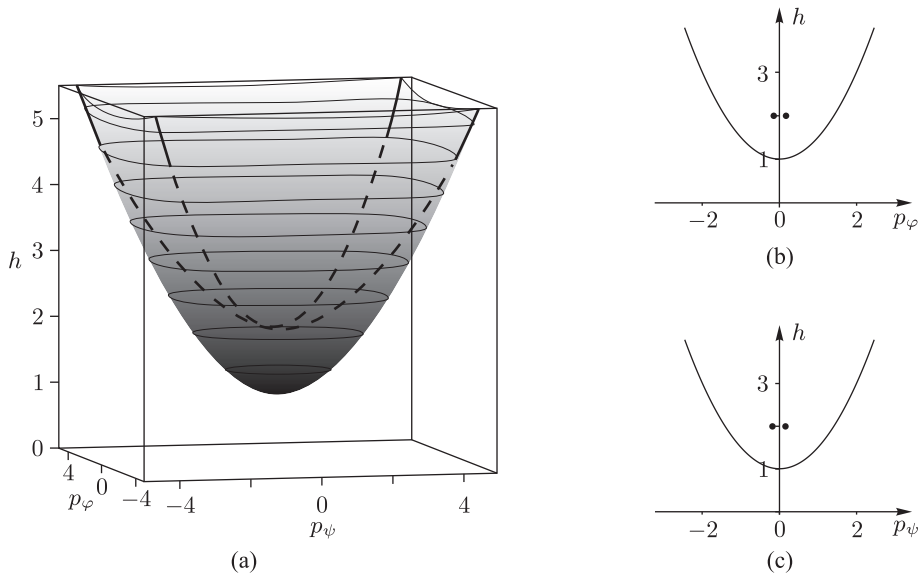}
\caption{Bifurcation diagram in the problem of an ellipsoid of revolution  on a smooth plane at $I_1=1$,
$I_3=1.5$, $b_1=1$, $b_2=2$, $m=1$, $g=1$:
(a)~three-dimensional bifurcation diagram in the~space $(p_{\varphi},p_{\psi}, h)$,
(b)~section formed by the intersection of the bifurcation diagram with the plane $p_{\psi}=0.157$,
(c)~section formed by the intersection of the bifurcation diagram with the
plane $p_{\varphi}=0.157$. }\label{fig_2}
\end{figure}

Fig.~\ref{fig_2} shows a three-dimensional bifurcation diagram of the system under consideration for the following parameters: $I_1=1$, $I_3=1.5$,
$b_1=1$, $b_2=2$, $m=1$, $g=1$. Visually the diagram consists of a ``bucket'' with two intersecting ``threads'' hanging above its bottom and corresponding to focus singularities.
At some critical value of energy~$h$ these ``threads'' reach the boundary of the ``bucket'' and become center-center singularities through the well-known Hopf bifurcation. From the viewpoint of dynamics
this phenomenon is well known as gyrostabilization of the top at high rotational velocities.
We can obtain two-dimensional bifurcation diagrams by intersecting the three-dimensional diagram with the plane $p_\psi=\const$ (or any other vertical plane). If the chosen plane
of section is not parallel to the planes $p_\psi=\pm p_\varphi$ and does not pass through the origin, then the corresponding two-dimensional bifurcation diagram contains
two focus singular points (see Fig.~\ref{fig_2}).

\subsection{Analysis of monodromy}\label{subsec_bkk1_3}

\begin{figure}[!ht]
\centering\includegraphics{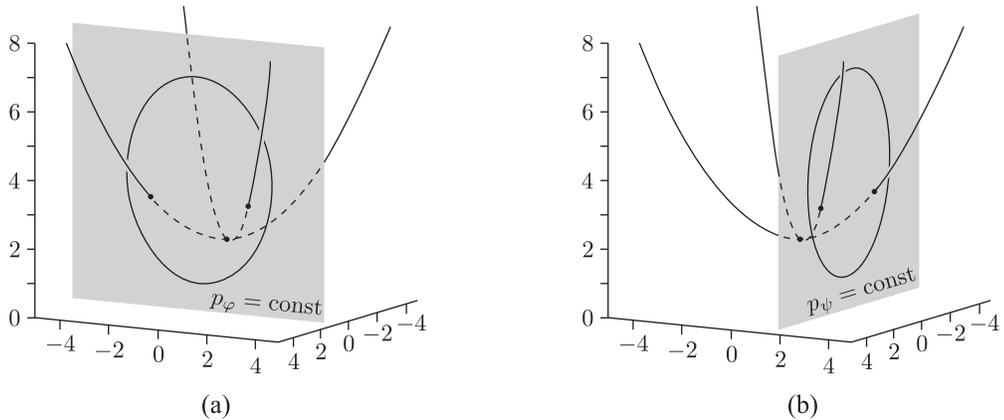}
\caption{Two scenarios for going around the focus singularities in the space $(p_{\varphi},p_{\psi},h)$.}\label{fig1*}
\end{figure}

To analyze the monodromy using the approach developed in Section~\ref{subsec2_1_},
it is necessary to choose a way of going around the singular points of the system in the space of first integrals $(p_{\varphi}, p_{\psi}, h)$. Whereas there is only one way
to go around one singularity, there exist two essentially different scenarios for going around two singularities.
One of them is a bypass at a fixed value of the integral $p_{\psi}$ (or a topologically equivalent bypass). The other
is a bypass at a fixed value of the integral $p_{\varphi}$ (or a topologically equivalent bypass). Fig.~\ref{fig1*} shows both bypass scenarios in the space $(p_{\varphi}, p_{\psi},
h)$. For a straightforward numerical analysis we choose the circles lying in the planes $p_{\psi}=\const$
and~$p_{\varphi}=\const$ as curves bypassing the singular points:
\begin{equation}
\label{1*}
\begin{aligned}
\gamma _{\psi}&= \big\{ (p_{\varphi}, p_{\psi}, h)\colon p_{\varphi}=p_{\varphi}^0+r_0\sin\alpha , \ \ p_{\psi}=p_{\psi}^0, \ \
    h=h^0+r_0\cos\alpha , \ \ \alpha \in[0,2\pi) \big\},\\
\gamma _{\varphi}&= \big\{ (p_{\varphi}, p_{\psi}, h)\colon p_{\varphi}=p_{\varphi}^0 , \ \ p_{\psi}=p_{\psi}^0 +r_0\sin\alpha, \ \
    h=h^0+r_0\cos\alpha , \ \ \alpha \in[0,2\pi) \big\}.
\end{aligned}
\end{equation}
The parameters $p_{\varphi}^0$, $p_{\psi}^0$, $h^0$ and~$r_0$ define the plane, the center of the circle and its radius and fix in the phase space the three-dimensional manifold $Q_3$
(see Section~\ref{subsec2_1_}), which is the preimage of the curve $\gamma _{\psi}$ (or $\gamma _{\varphi}$) for integral mapping. The variable $\alpha $ is an angle coordinate
on the curve $\gamma _{\psi}$ (or $\gamma _{\varphi}$) and parameterizes
the family of invariant tori forming $Q_3$. On each of the invariant tori it is convenient to choose as coordinates the Euler angles $\varphi$ and~$\theta$, which are convenient for
calculations, although they are not classical angle coordinates on a torus.
A schematic of the manifold is shown in Fig.~\ref{fig3_}.

On the invariant torus corresponding to some value $\alpha$,  as basis cycle $\mu_{\alpha }$ (invariant under the symmetry field of the system) we choose a cycle
given by the relation
$\dot\theta=\const\ne0$.

The family of cycles $\mu_{\alpha }$, $\alpha =[0,2\pi)$ forms the torus
$$
T_{\rm transv}=\{{\bs x}\colon {\bs x}\in\mcM^5, \Phi({\bs x})\in \gamma _{\psi},\dot\theta=\const\}.
$$
As stated above, the phase flow of the system is transverse to this torus and defines the Poincar\'{e} map~\eqref{Poincaremap} on it. As angle coordinates on the torus $T_{\rm transv}$
one can choose  $\alpha $ and~$\varphi$.
The basis cycle $\mu_{\alpha }$ in these coordinates is a vertical straight line $\alpha =\const$ and is obviously invariant under the Poincar\'{e} map, since the value $\alpha $
is the integral of motion of the system under consideration. Thus, as the first basis cycle on the torus $T_{\rm transv}$ we can choose the cycle $\alpha =\const$.

We now define the second basis cycle $\nu$ on the torus $T_{\rm transv}$ by $\varphi=0$. Thus, the analysis of monodromy reduces to analysis of iterations of the cycle $\nu=\{\varphi=0\}$
with the Poincar\'{e} map~\eqref{Poincaremap} on the torus $T_{\rm transv}=\{(\alpha
,\varphi)\}$. The nontrivial monodromy with coefficient $k$
(see Section~\ref{subsec2_1_}) corresponds to the case where the image $\sigma(\nu)$ of the basis cycle
$\nu$ is a curve making $k$ turns in the direction $\varphi$ as~$\alpha $ changes from zero to $2\pi$.

\begin{figure}[!ht]
\centering\includegraphics{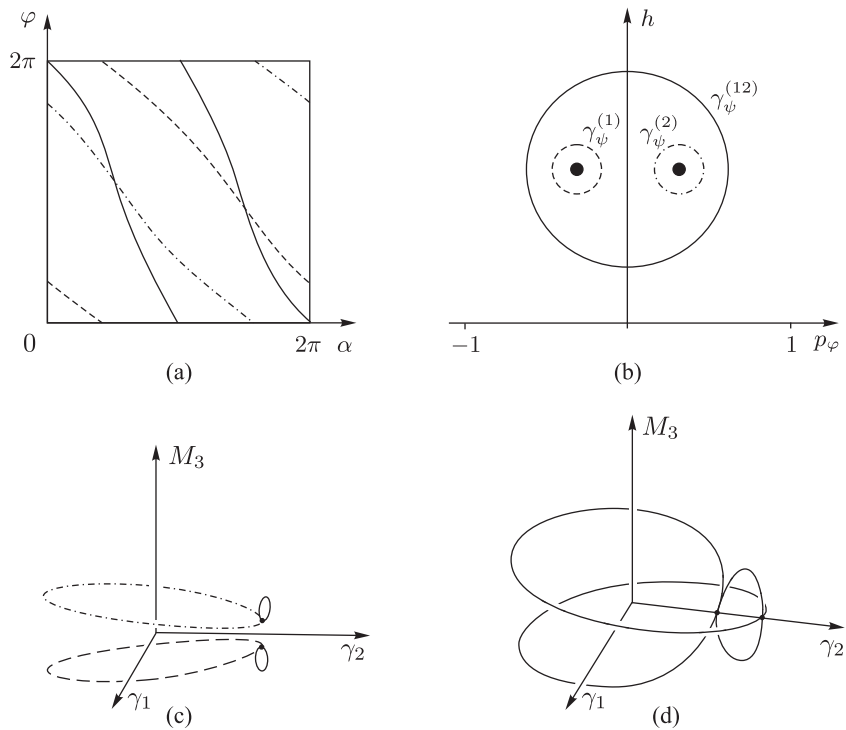} \caption{Images $\sigma(\nu)$ of the basis cycle $\nu=\{\varphi=0\}$ under the Poincar\'{e} map:
(a)~on the torus $T_{\rm transv}$,
(c,\,d)~in the space $(\gamma _1,\gamma _2,M_3)$;
(b)~the corresponding curves going around the singularities in the~plane $p_{\psi}=0.157$.}\label{fig2*}
\vspace{-2mm}
\end{figure}

We now consider the results of numerical investigation of the Poincar\'{e} map~\eqref{Poincaremap} for the system~\eqref{eq1_1}. Fig.~\ref{fig2*}a shows the
torus $T_{\rm transv}$ as a square in coordinates $\alpha $ and~$\varphi$. It can be seen that the basis cycle $\nu$
is the straight line $\varphi=0$. The curves depicted on the square are the images of this basis cycle under the action of the Poincar\'{e} map~\eqref{Poincaremap} in cases
of going around the singularities of the system along three different curves $\gamma
_{\psi}^{(1)}$, $\gamma _{\psi}^{(2)}$, $\gamma _{\psi}^{(12)}$ lying in the
plane $p_{\psi}=\const$ and depicted in Fig.~\ref{fig2*}b. The curves $\gamma _{\psi}^{(1)}$ and~$\gamma _{\psi}^{(2)}$ correspond to a ``walk'' around
each single singularity, and~$\gamma _{\psi}^{(12)}$~correspond to simultaneously going around both singularities.
The same images of basis cycles in the space $(\gamma _1,\gamma _2,M_3)$
are presented in Fig.~\ref{fig2*}с (for the curves $\gamma _{\psi}^{(1)}$ and~$\gamma
_{\psi}^{(2)}$) and in Fig.~\ref{fig2*}d (for the curve $\gamma _{\psi}^{(12)}$).
The figure shows that in the case of going around each of the foci the monodromy is nontrivial ($k=-1$), while in the case of simultaneously going around both foci it doubles.

\begin{figure}[!ht]
\centering\includegraphics{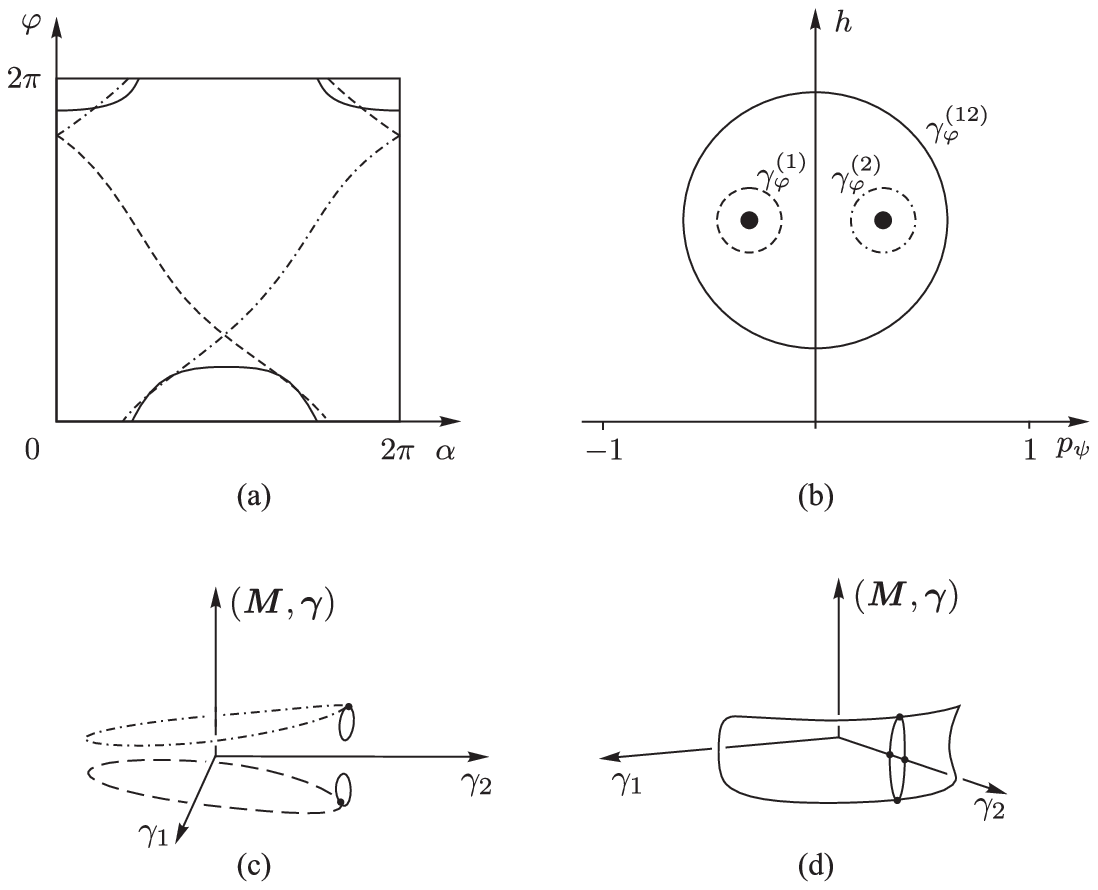} \caption{Images $\sigma(\nu)$ of the basis cycle $\nu=\{\varphi=0\}$
for the Poincar\'{e} map: (a)~on the torus $T_{\rm transv}$,
(c,\,d)~in the space $\big(\gamma _1,\gamma _2,({\bs M},\bgam)\big)$;
(b)~the corresponding curves going around the singularities in the~plane $p_{\varphi}=0.157$.}
\label{fig3*}
\vspace{-3mm}
\end{figure}

Fig.~\ref{fig3*} presents analogous calculation results for the curves $\gamma _{\varphi}^{(1)}$, $\gamma _{\varphi}^{(2)}$
and~$\gamma_{\varphi}^{(12)}$ lying in the~plane $p_{\varphi}=\const$. As seen in the figure, the curve $\gamma _{\varphi}^{(1)}$ corresponds to the negative monodromy $k=-1$,
the curve $\gamma _{\varphi}^{(2)}$~to the positive monodromy $k=+1$,
and for the curve~$\gamma _{\varphi}^{(12)}$ the monodromy is trivial.

\begin{figure}[!ht]
\centering\includegraphics{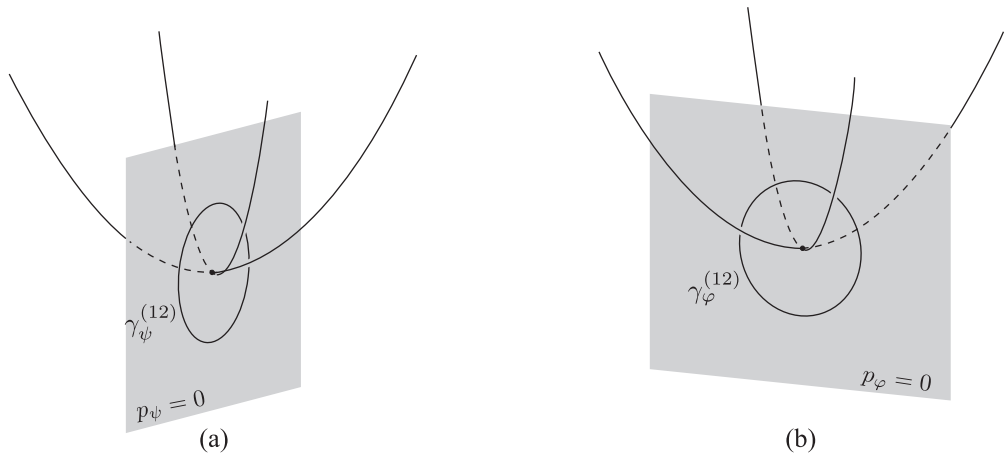}
\caption{Two ways to go around a double focus point.}
\label{figure_3*}
\end{figure}

The cases considered here correspond to a walk around two focus singularities at small but nonzero values of the integrals $p_\varphi$ and~$p_\psi$. If their values are assumed
to be equal to zero, the pairs of focus points shown in Figs.~\ref{fig2*}b and~\ref{fig3*}b merge into a single one. In the phase space $\mathcal M^5$ this corresponds
to the situation where the focus points are at the same singular level of integral mapping,
that is, this level becomes a double pinched torus. On the three-dimensional bifurcation diagram this event corresponds to the point of intersection of two hanging
threads.  From the viewpoint of monodromy, of course, there will be no changes, since the circles $\gamma_\psi^{(12)}$ and
$\gamma_\varphi^{(12)}$ will undergo no bifurcations (see Fig.~\ref{figure_3*})

In other words, when going around a torus with two pinch points in the~plane $p_\psi =
0$, we obtain a monodromy matrix of the form $\begin{pmatrix}  1 & 2 \\[-2mm] 0 & 1
\end{pmatrix}$, as it should be in the Hamiltonian case, and  when a point analogously goes around the same fiber in the~plane
$p_\varphi = 0$ the monodromy matrix turns out to be trivial, i.e.,
$\begin{pmatrix}  1 & 0 \\[-2mm] 0 & 1 \end{pmatrix}$, which is impossible in the Hamiltonian case (see Section~\ref{sec_bkk2}). Of course, there is no
contradiction here. The statement that in the Hamiltonian case the monodromy matrix must have the form
$\begin{pmatrix}  1 & 2 \\[-2mm] 0 & 1
\end{pmatrix}$ applies to Hamiltonian systems defined on four-dimensional
symplectic manifolds, i.e., a walk around a singular fiber must be performed on the  {\it symplectic} level.
The level $p_\psi = 0$ satisfies this condition, as it is a symplectic fiber of the Poisson structure. On the contrary,
the level $p_\varphi=const$
possesses no symplectic structure, and therefore no Hamiltonian monodromy occurs in this case.

\section{Rolling of an ellipsoid of revolution on a rough plane}\label{sec4}

\subsection{Equations of motion and first integrals}\label{subsec2_1}

Now consider an axisymmetric ellipsoid rolling on a plane without slipping under the same assumptions about its form and mass distribution as in Section~\ref{sec_bkk1}.
The equations of motion, first integrals and a discussion of the dynamics of an ellipsoid (and an arbitrary body) rolling on a plane can be found in~\cite{BorMamRolling}.
Also notice that the bifurcation analysis for quite a similar problem of a flat round disk rolling on a plane, which can be treated as a limiting case of an oblate ellipsoid of revolution, was done in \cite{20}.

The equations of motion in our case are
\begin{equation}
\label{eq12}
\left\{
\begin{aligned}
&\dot{{\bs M}}={\bs M}\times\bom-m\dot{{\bs r}}\times(\bom\times{\bs r})+mg{\bs r}\times\bgam,\\
&\dot{\bgam}=\bgam\times\bom,
\end{aligned}
\right.
\end{equation}
where the angular momentum of the body~${\bs M}$ w.r.t. the point of contact is related to the angular velocity of the ellipsoid~$\bom$ by
\begin{equation}
\label{eq13}
{\bs M}=\bfI\bom+m{\bs r}\times(\bom\times{\bs r}),
\end{equation}
and~${\bs r}$ and~$\bgam$ are expressed in terms of each other with the help of~\eqref{eq1_3} as before.

The equations~\eqref{eq12} admit the energy integral
\begin{equation}
\label{eq14}
H=\dfrac{1}{2}\,({\bs M},\bom)-mg({\bs r},\bgam)
\end{equation}
and the geometric integral~$\bgam^2=1$, and define, as in the previous case, the flow on~$\mcM^5$. In addition, in the axisymmetric case the equations of motion
admit an invariant measure of the form
\begin{equation}
\label{eq15}
\rho=\dfrac{1}{\sqrt{I_1I_3+m({\bs r},\bfI{\bs r})}}
\end{equation}
and two more integrals of motion. In the case of an arbitrary body of revolution these integrals are linear in~${\bs M}$, however, they are expressed in terms of non-algebraic
(for a disk, for example,~hypergeometric) functions of~$\bgam$.
Below we present a short algorithm for obtaining these integrals (for details see~~\cite{BorMamRolling}).

We introduce the variables~$K_1$ and~$K_2$ invariant under rotations about the symmetry axis
\begin{equation}
\label{eq16}
\begin{aligned}
K_1&=M_1\gamma _1+M_2\gamma _2+\dfrac{b_3^2}{b_1^2}\,M_3\gamma _3,\\
K_2&=\dfrac{\omega _3}{\rho}=\rho\left(\dfrac{mb_1^2b_3^2\gamma _3}{b^{2}_{1}+
(b^{2}_{3}-b^{2}_{1})\gamma ^{2}_{3}}\left(M_1\gamma _1+M_2\gamma _2+
\dfrac{b_3^2}{b_1^2}\,M_3\gamma _3\right)+I_1M_3\right)\!\!,
\end{aligned}
\end{equation}
and as the new time we choose~$\gamma _3$. Then the equations of motion for~$K_1$, $K_2$ take the form of the linear system
\begin{equation}
\label{eq17}
\begin{gathered}
K'_1=\rho I_3\dfrac{b_3^2-b_1^2}{b_1^2}K_2,\quad
K'_2=\dfrac{m\rho b_1^4(b_3^2-b_1^2)(1-\gamma _3^2)}{\big(b_1^2+(b_3^2-b_1^2)\gamma _3^2\big)^2}K_1.
\end{gathered}
\end{equation}
The general solution to this system can be represented as
\begin{equation}
\label{eq18}
\bs K =\bfG(\gamma _3)\bs C,
\end{equation}
where $\bs K =(K_1,K_2)$, $\bfG$~is the fundamental solution matrix of the system~\eqref{eq17} with the initial conditions $\bfG(\gamma _3=0)=\mathrm{Id}$, which can be expressed
in terms of the Heun functions, and $\bs C=(C_1,C_2)$~are the constants of integration, which are the sought-for first integrals.
As above, by the capital letters~$C_1$ and~$C_2$  we denote the integrals of motion as functions of the phase variables and by the lowercase letters~$c_1$ and~$c_2$~the
values which they take on specific trajectories.
By virtue of the uniqueness theorem the matrix~$\bfG(\gamma _3)$ is reversible for all values of~$\gamma _3$. Thus, the first integrals of motion have the form
\[
\bs C=\bfG^{-1}(\gamma _3)\bs K
\]
and are functions linear in momenta, with coefficients which are non-algebraic functions of~$\gamma _3$.
We note that the integrals~$C_1$ and~$C_2$ are not equivalent.
Analysis of the expressions~\eqref{eq16} for determination of~$K_1$ and~$K_2$, and of the chosen initial conditions of the fundamental solution matrix
$\bfG(\gamma _3=0)=\mathrm{Id}$ suggests that the integral~$C_1$ is an analog of the area integral  $F_1$ in the previous problem and the integral~$C_2$~is an analog of the Lagrange integral $F_3$.
In what follows we shall ascertain the correctness of this analogy when analyzing the monodromy of the systems under consideration.

\subsection{Bifurcation analysis}

As in the previous section, we introduce the local coordinates ${\bs x}=(\gamma _3, \dot\gamma _3, \varphi,C_1,C_2)$. By calculating the corresponding
Jacobian $\frac{\partial {\bs x}}{\partial ({\bs M},\bgam)}$, it is easy to show that these coordinates are defined everywhere on $\mcM^5$
except for  those points where $\gamma _3=\pm1$. Using the approach developed in Section~\ref{subsec_bkk1_2}, it is straightforward to show that
the critical set $S_2$ is given by two equations
\begin{equation}
\label{eq19}
\dfrac{\partial H}{\partial \gamma _3}=0,\quad
\dfrac{\partial H}{\partial \dot\gamma _3}=0,
\end{equation}
where $H$~is the energy integral expressed in terms of the chosen coordinates and has the form
\begin{equation}
\label{eq20}
H=\dfrac1{2I_1(1-\gamma _3^2)}\left(K_1^2-\dfrac{I_3l^2K_2^2}{mb_1^2}K_2^2 +\dfrac{mb_3^2\gamma_3^2}{I_1l^2}
    \bigg(K_1-\dfrac{l^2K_2}{\rho mb_1b_3}\bigg)^2 \right)+
    \dfrac12 k^2\dot\gamma _3^2 +U(\gamma _3^2),
\end{equation}
$\displaystyle l=\sqrt{b_1^2(1-\gamma_3^2)+b_3^2\gamma_3^2}$~is the height of the center of mass of the ellipsoid, and~$K_1$ and~$K_2$  are expressed in terms of the
integrals $C_1$,
$C_2$ and~$\gamma _3$ using~\eqref{eq18}. It follows from the second of Eqs.~\eqref{eq19} that $\dot\gamma _3=0$. Thus, the bifurcation surface in this case is given by the
relation
\begin{equation}
\label{eq21}
H=\dfrac1{2I_1(1-\gamma _3^2)}\left(K_1^2-\dfrac{I_3l^2K_2^2}{mb_1^2}K_2^2 +\dfrac{mb_3^2\gamma_3^2}{I_1l^2}
    \bigg(K_1-\dfrac{l^2K_2}{\rho mb_1b_3}\bigg)^2 \right)+U(\gamma _3^2),
\end{equation}
where $\gamma _3$ is the solution of the equation $\frac{\partial H}{\partial \gamma _3}\Big|_{\dot\gamma _3=0}=0$.

As in the case of the ellipsoid rolling on a smooth plane, the case $\gamma
_3=\pm1$ should be considered separately.
Obviously there exist two types of trajectories for which the equality $\gamma _3=\pm1$ holds. The first type includes trajectories transversally intersecting the submanifold
$\{{\bs x}\colon \gamma _3=\pm1\}$. For such trajectories the relations $\dot\gamma _3\big|_{\gamma _3=\pm1}=0$, $\ddot\gamma _3\big|_{\gamma _3=\pm1}\ne 0$ are satisfied.
Due to continuity, for the motion along such trajectories in an arbitrarily small neighborhood of point $\gamma _3=\pm1$ the derivative $\dot\gamma _3$
becomes different from zero. Consequently, such a trajectory cannot be a singular periodic solution, since it requires that the equality $\dot\gamma _3=0$ be satisfied.

\begin{figure}[!ht]
\centering\includegraphics{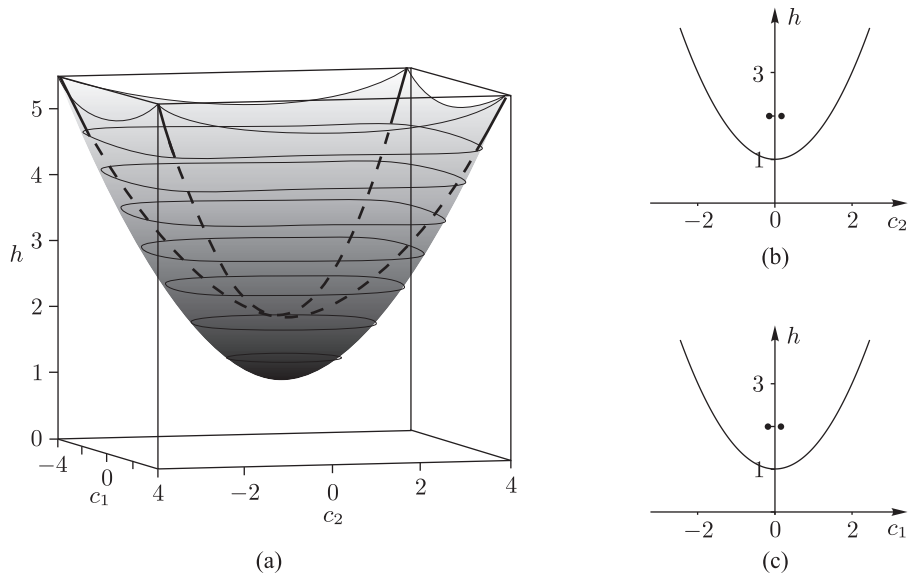}
\caption{Bifurcation diagram in the problem of an axisymmetric ellipsoid rolling on a rough plane for $I_1=1$,
$I_3=1.5$, $b_1=1$, $b_2=2$, $m=1$, $g=1$:
(a)~three-dimensional bifurcation diagram in the space $(c_2,c_1, h)$,
(b)~section formed by the intersection of the bifurcation diagram with the plane $c_1=0.157$,
(c)~section formed by the intersection of the bifurcation diagram with the plane $c_2=0.157$.}\label{fig3}
\end{figure}

The second type of trajectories for which the equality $\gamma
_3=\pm1$ holds includes trajectories completely lying on this submanifold. For these trajectories the relations $\gamma_3=\pm1$, $\dot\gamma_3=0$ and $\ddot\gamma_3=0$ hold.
Differentiating the second equation~\eqref{eq12} and substituting the equalities $\gamma _1=\gamma_2=0$ into it, we obtain
$$
\ddot\gamma_3\Big|_{\gamma _3=\pm1}=-\frac{M_1^2+M_2^2}{(I_1+mb_3^2)^2},
$$
hence, for solutions of the second type we have
${\bs M}=(0,0,M_3)$ and such solutions are vertical rotations about the symmetry axis.
It is straightforward to show that the corresponding Jacobi matrix has the form
$$
\dfrac{\partial(H,C_1,C_2)}{\partial ({\bs M},\gamma _1,\gamma _2)}\Bigg|_{\substack{\gamma _3=\pm1\\{\bs M}=(0,0,M_3)}}=
\begin{pmatrix}
0 & 0 & \frac{\partial H}{\partial M_3} & 0 & 0\\
0 & 0 & \frac{\partial C_1}{\partial M_3} & 0 & 0\\
0 & 0 & \frac{\partial C_2}{\partial M_3} & 0 & 0
\end{pmatrix}\!\!,
$$
and its rank drops by two. Thus, as in the case of a smooth ellipsoid, the vertical rotations about the symmetry axis form two one-parameter families of singular points
$S_1 =\{{\bs M}=(0,0,M_3),
\bgam=(0,0,\pm1)\colon M_3\in(-\infty,\infty)\}$, whose image in the space of first integrals consists of two curves. The corresponding bifurcation diagram
is presented in Fig.~\ref{fig3} and does not qualitatively differ from the diagram for the case of a smooth plane.

For a generic ellipsoid with distinct semiaxes  (even for a homogeneous one), the situation is much more complicated. So far, it is not even clear whether there exists an invariant measure
in this case. It is known that for completely asymmetric bodies, its absence leads to the existence of a strange attractor and contradicts the property of being Hamiltonian and
conformally Hamiltonian \cite{27,28,29}.

\subsection{Analysis of monodromy}

\begin{figure}[!ht]
\centering\includegraphics{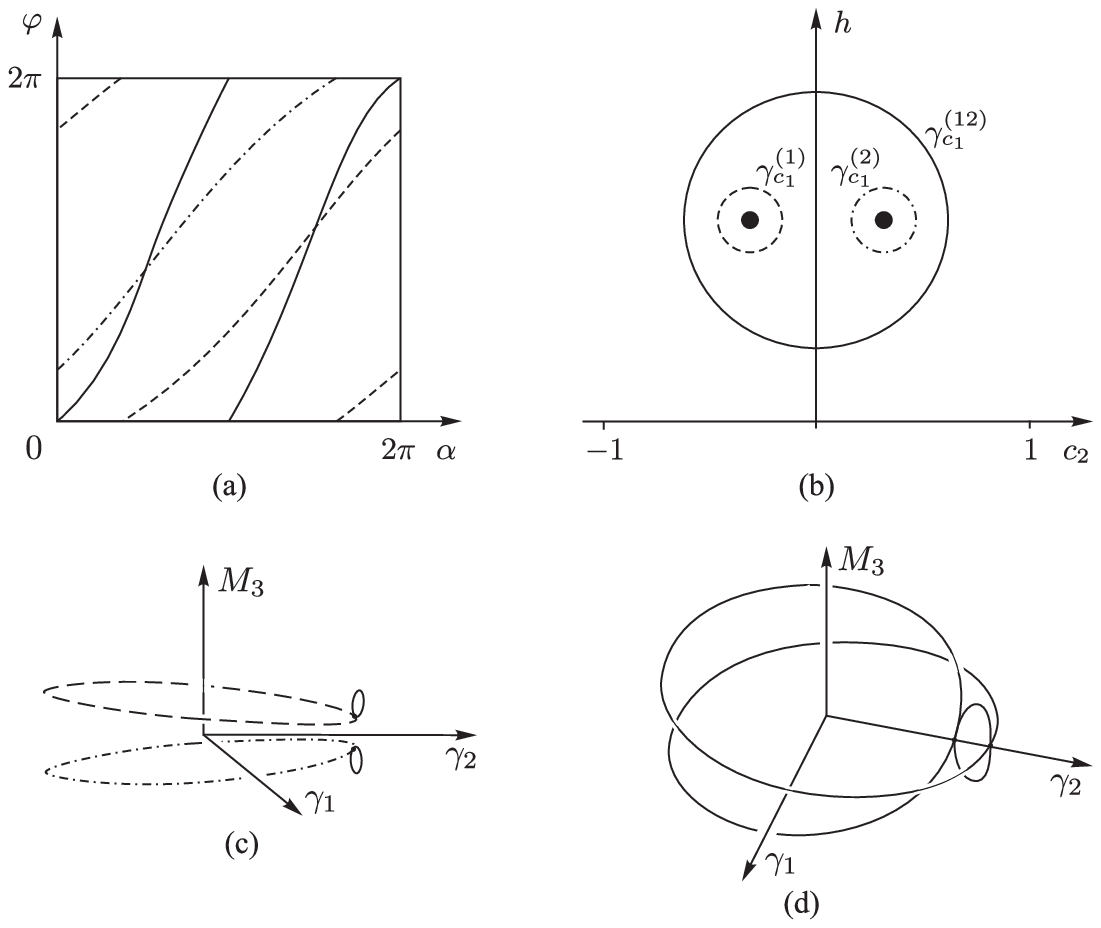} \caption{
Images $\sigma(\nu)$ of the basis cycle $\nu=\{\varphi=0\}$ under the Poincar\'{e} map: (a)~on the torus $T_{\rm transv}$,
(c,\,d)~in the space $\big(\gamma _1,\gamma _2,M_3\big)$;
(b)~the corresponding curves going around the singularities in the~plane $p_{\varphi}=0.157$.}\label{fig7}
\end{figure}

\begin{figure}[!ht]
\centering\includegraphics{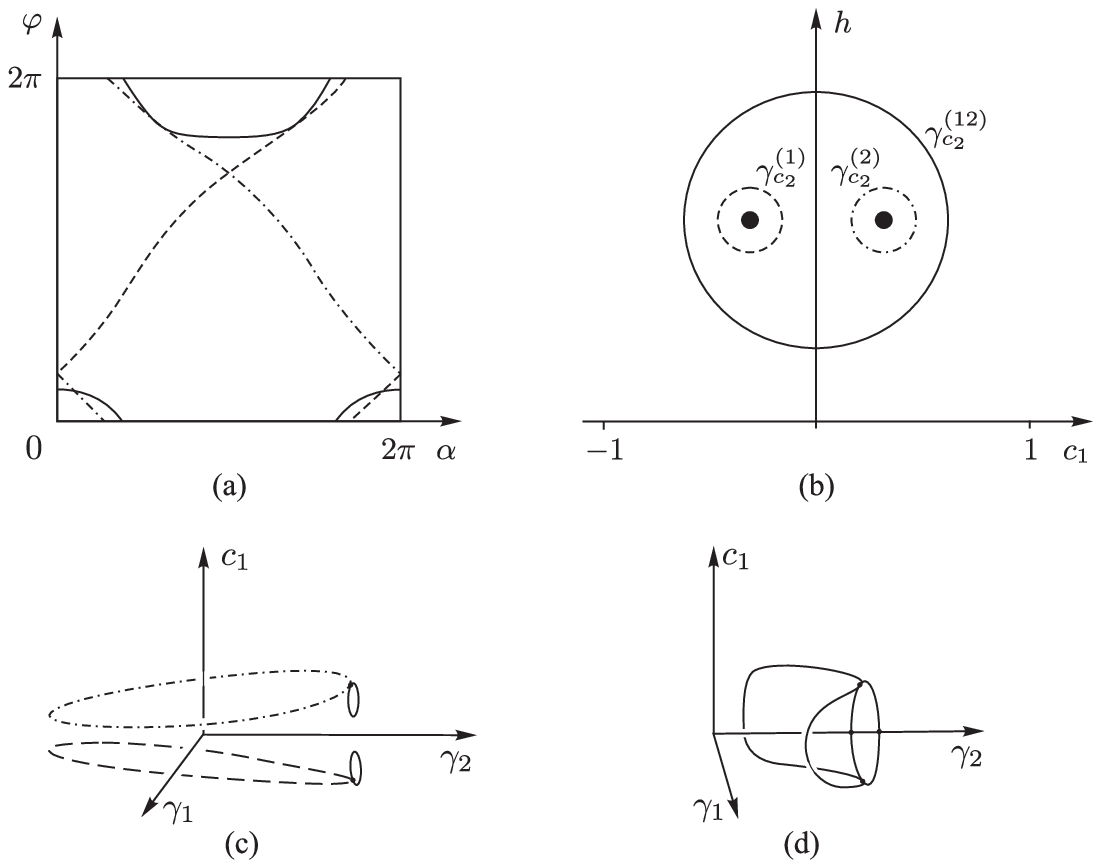}
\caption{
Images $\sigma(\nu)$ of the basis cycle $\nu=\{\varphi=0\}$
under the Poincar\'{e} map: (a)~on the torus $T_{\rm transv}$,
(c,\,d)~for the Poincar\'{e} map in the space $\big(\gamma _1,\gamma _2,c_1\big)$;
(b)~the corresponding curves going around the singularities in the~plane $p_{\varphi}=0.157$.}\label{fig8}
\end{figure}

To analyze the monodromy, we use the approach developed in Section~\ref{subsec_bkk1_3}. It is sufficient to replace $p_{\psi}$
with $c_1$ and ~$p_{\varphi}$ with $c_2$ in our reasoning. The results of building a map for the
curves lying in the~plane $c_1=\const$ are presented in Fig.~\ref{fig7} and those for the curves lying in the~plane $c_2=\const$ are shown in Fig.~\ref{fig8}. As seen in the
figures, the monodromy for the curves in the~plane $c_1=\const$ is identical (up to sign) with the case $p_{\psi}=\const$, and for the curves lying in the plane $c_2=\const$ with
the case $p_{\varphi}=\const$ for the ellipsoid rolling on a smooth plane.

Thus, despite the fact that the  problem of an ellipsoid of revolution rolling on a rough plane is nonholonomic, its monodromy completely coincides with the one of
the Hamiltonian problem of rolling on a smooth plane. Moreover, these two problems
{\it do not differ from each other topologically} at all, which allows us to conjecture that the first (nonholonomic) problem is, in fact,  conformally Hamiltonian.

\section{Results of the analysis and conclusions}

The main result of our analysis is the confirmation of our counter conjecture:
the nonholonomic integrable system describing the dynamics of an ellipsoid of revolution on a {\it rough} plane, in its topological properties, is quite analogous to the Hamiltonian
system describing the dynamics of the same ellipsoid on a {\it smooth} plane.  In other words, the monodromy gives no topological obstructions to the Hamiltonization of the
nonholonomic system in question.

However, this conclusion does not mean that the monodromy is useless for the Hamiltonization problem.  On the contrary, the above example demonstrates its
exceptional effectiveness. Indeed,  when solving the Hamiltonization problem we usually want to construct a Poisson (but not a symplectic!) structure
relative to which a system under consideration turns out to be Hamiltonian. It is natural to require the energy integral known a priori to be the Hamiltonian of the system. The situation
with the other integrals is not so clear. One needs to ``partition'' them somehow into Casimir functions and ``real'' integrals. The choice of the Casimir functions
in this context is equivalent to defining the foliation of the space $\mathcal M^5$ into symplectic leaves.

The monodromy helps to eliminate an incorrect ``partition''.  It is this phenomenon that we observe in the problem under consideration. In addition to the energy integral, our nonholonomic system
possesses two linear integrals, which at first glance do not considerably differ from each other. For Hamiltonization of the system, it is natural to ``choose'' one of them as
a Casimir function of the sought-for Poisson structure.  Which of the two? The analysis shows that $C_2$ does not suit for this purpose,  since with this choice the monodromy
around the singular fiber becomes non-Hamiltonian.
Conversely,  $C_1$ is quite suitable: from the topological point of view the foliation into hypothetic symplectic leaves will be ``like'' the standard one, and there will be no problems with
monodromy.\footnote{Non-equivalence of the integrals can be detected in another way as well. The levels of the first of them $\{
C_1=\mathrm{const}\}$ are diffeomorphic to the cotangent bundle to the sphere
$T^*S^2$,  whereas for the other integral $C_2$ they will be direct products of $S^2\times \R^2$ (i.e., trivial $\R^2$-bundles over the sphere).  However,  this does not cause any problem from the symplectic point of view: both are good symplectic manifolds.}

By the way, an explicit Hamiltonization of an ellipsoid of revolution on a
rough plane is still not accomplished (we mean the existence of a suitable Poisson structure of rank~4; the structure of rank~2 for this problem after some additional reduction
has been found in~\cite{BorMamRolling}), although
attempts to find an explicit conformally Hamiltonian representation for this kind of problems have been made.  However,  even for a simpler problem of Routh's sphere (a dynamically symmetric ball with a displaced center of mass), so far one has not succeeded in constructing a Poisson structure without singularities~\cite{23}.

We also note that the method for calculating the monodromy using the Poincar\'{e} map has proved to be very effective, illustrative and stable. This method can also be applied
successfully in studying other (even non-integrable) problems, since the Poincar\'{e} map can be defined for much more general dynamical systems.

\bigskip

The authors thank A.\,V.\,Borisov, I.\,S.\,Mamaev and~B.\,Zhilinskii for useful discussions. A.\,V.\,Bolsinov would like to emphasize an exceptional importance
of the discussions with H.\,Dullin on topological and dynamical aspects of the monodromy phenomenon  in June 2009,  many ideas from those conversations  have been used in this work
in the most immediate way.

\bigskip

This research was supported by the Federal Target Program ``Scientific and Scientific-Pedagogical Personnel of Innovative Russia'' for 2009--2013,  Agreement  No\,14.В37.21.1935
``Topological Methods in Mechanics  and Hydrodynamics'', 	the Analytical Departmental Target Program ``Development of Scientific Potential of Higher Schools`` for 2012--2014, No\,1.1248.2011
``Nonholonomic Dynamical Systems and Control Problems'', the Analytical Departmental Target Program ``Development of Scientific Potential of Higher Schools'' for 2012--2014,
No\,1.7734.2013
``Development of Nonholonomic Mobile Systems''.
The work of A.\,A.\,Kilin was supported by the grant of the President of the Russian Federation for young Doctors of Science MD-2324.2013.1
``Nonholonomic Mobile Systems: Various Models of Friction, Regular and Chaotic Regimes and Control''.

\mbox{}\label{lpage}
\end{document}